\documentclass[12pt]{amsart}

\usepackage{amssymb,amsfonts,amsmath,color}

\begingroup
\newtheorem{thm}{Theorem}

\begin{document}

\title[String $C$-groups with real Schur index $2$]{String $C$-groups with real Schur index $2$} 

\author{Peter J.~Cameron} 
\address{
School of Mathematics and Statistics, 
University of St Andrews,
North Haugh,
St Andrews, Fife KY16 9SS,
Scotland
}
\email{pjc20@st-andrews.ac.uk} 

\author{Allen Herman$^*$}\thanks{$^*$The work of this author is supported by an NSERC Discovery Grant.}
\address{Department of Mathematics and Statistics, University of Regina, Regina, Canada, S4S 0A2}
\email{Allen.Herman@uregina.ca}

\author{Dimitri Leemans$^\circ$}\thanks{$^\circ$The work of this author is supported by an ARC Advanced grant of the Communaut\'e Fran\c caise Wallonie-Bruxelles.}
\address{Universit\'e Libre de Bruxelles, D\'epartement de Math\'ematique, C.P.216 - Alg\`ebre et Combinatoire, Boulevard du Triomphe, 1050 Brussels, Belgium}
\email{dleemans@ulb.ac.be} 



\begin{abstract}
We give examples of finite string $C$-groups (the automorphism groups of abstract regular polytopes) that have irreducible characters of real Schur index $2$.  This answers a problem of Monson concerning these groups.
\end{abstract}

\date{April 13, 2021}
\maketitle

\noindent {\bf 1.  Monson's question.} 

\bigskip
In 1990 Barry Monson observed that some counting formulas obtained by Peter McMullen involving irreducible representations of string $C$-groups (aka.~the groups of automorphisms of finite abstract regular polytopes) did not take into account the possibility that irreducible constituents of the character afforded by the geometric simplex realization of the abstract regular polytope might not be realizable over the real numbers.  Both situations of this, either an imaginary field of character values, or a real field of character values but not realizable over the reals, were not considered by McMullen.  A real-valued irreducible character that is not realizable over the field of real numbers is one with real Schur index $2$.  The formulas were revised to account for this in \cite{McMullen-Monson}, but existence of the characters in question was not resolved there.  The existence of complex-valued irreducible characters was observed in \cite{Herman-Monson}. Some of the earlier misperceptions of the interplay between representation theory and geometry for string $C$-groups persisted until very recently clarified by Frieder Ladisch \cite{Ladisch15}, but existence of irreducible characters of real Schur index $2$ for these groups remained unanswered.  

Monson's question appeared as Problem 23 in Schulte-Weiss' {\it Problems on Polytopes} \cite{SW06}.  

\medskip
\noindent {\bf Problem 23:} Can one finite string $C$-group have an irreducible character of real Schur index $2$? 

\medskip
The important irreducible characters of string $C$-groups for discrete geometry are the {\it polytopal} ones.  These are the irreducible constituents (other than $1_G$) of the permutation character $(1_H)^G$ where $H$ is the stabilizer of a vertex in the regular abstract polytope whose automorphism group is $G$, and they are precisely the  irreducible constituents of the simplex realization of the polytope with automorphism group $G$ for the given string $C$-group representation.  Polytopal irreducible characters are the only ones that occur in McMullen's formulas for the dimension of the realization cone of the polytope.  So if the answer to Problem 23 is yes then it makes sense to ask:   

\medskip
\noindent {\bf Question:} If there are irreducible characters of a string $C$-group with real Schur index $2$, can they occur as polytopal irreducible characters with respect to some string $C$-group representation? 

\medskip
A recent article of Herman and Alsairafi explains why all polytopal irreducible characters are of real type for all finite string Coxeter groups and for all string $C$-groups that are quotients of rank $3$ affine Coxeter groups of type $\{4,4\}$ and $\{3,6\}$ \cite{Herman-Alsairafi}.  In this article, we show there are infinitely many projective special unitary groups $PSU(n,q)$ (aka.~the simple groups $U_n(q)$ or the simple groups of Lie type ${}^2A_{n-1}(q)$), that have string $C$-group representations having some irreducible characters of real Schur index $2$, and for some of these groups all of their irreducible characters are polytopal relative to at least one string $C$-group representation.   

The smallest such group is $PSU(5,3)$, for this group both of these facts can be verified computationally using the information in GAP's Atlas Library \cite{GAP}.  In Section 3 we outline direct GAP calculations showing $PSU(5,3)$ has a string $C$-group redrepresentation of Schl\''{a}fi type $\{5,5\}$ relative to which every irreducible character of the group is polytopal.  So this group answers both questions affirmatively. 

That such simple groups of Lie type are string $C$-groups (of rank $3$) is a result of Nuzhin \cite{Nuzhin}.  That they have irreducible characters of real Schur index $2$ was observed by Ohmori \cite{Ohmori}, who showed this is the case for their cuspidal unipotent irreducible characters.  

\bigskip
\noindent {\bf 2. String $C$-groups.}

\medskip
Let $\Gamma$ be an infinite Coxeter group of rank $r=n+1$, with Schl\"afli type $\{m_1,\dots,m_n\}$.  This implies that $\Gamma$ has involutory generators $\rho_0,\rho_1,\dots,\rho_n$ whose relations are $(\rho_{i-1}\rho_i)^{m_i} = 1$ and $(\rho_j\rho_i)^2=1$ if $|i - j|\ge 2$.  A {\em string $C$-group of Schl\"afli type $\{m_1,\dots,m_n\}$} is a finite quotient $G = \Gamma/N$ where $N$ is a normal subgroup of finite index for which the set of involutary generators $S = \{s_0, s_1, \dots, s_n : s_i :=\rho_i N \mbox{ for } i = 0,1,\dots,n\}$ satisfies the intersection condition: for every pair of subsets $I,J \subseteq S$, $\langle I \rangle \cap \langle J \rangle = \langle I \cap J \rangle$.   The pair $(G,S)$ is then called a {\em string C-group representation} of $G$ or {\em string C-group} for short and when the context is clear the set $S$ is omitted. The {\em rank} of $(G,S)$ is the cardinality of $S$. 

The string $C$-group $G$ is the automorphism group of an abstract regular polytope $\mathcal{P}$. The automorphism action is realized by a transitive action on the collection of maximal-length flags of the polytope.  The stabilizer of a vertex in this action is a conjugate of the subgroup $H = \langle s_1,\dots,s_r \rangle$.  The simplex realization can be interpreted as in \cite{Herman-Monson} to be a representation of the group $G$ whose character affords the deleted permutation character $\sigma = (1_H)^G-1_G$.  The irreducible constituents of $\sigma$ are precisely the polytopal irreducible characters of the string $C$-group $G$ with respect to the defining string $C$-group representation.   

\bigskip
Problem 23 asks about the existence of (polytopal) irreducible characters of real Schur index $2$.  For a specific irreducible character, this property can be checked directly from the calculation of the Frobenius-Schur indicator $$\nu_2(\chi) := \frac{1}{|G|} \sum_g \chi(g^2),$$ 
since $m_{\mathbb{R}}(\chi)=2 \iff \nu_2(\chi) = -1$ (see \cite{Isaacs}).   

\bigskip
\noindent {\bf 3. Calculations for $PSU(5,3)$.}

\medskip
In this section we show how to directly verify our claims using the Atlas \cite{Atlas} via GAP's {\tt AtlasRep} package \cite{AtlasRep}.  For standard GAP functions we do not give details but refer readers to the GAP manual \cite{GAP}.  For {\sc Magma}~\cite{BCP97} users, there is a similar Atlas interface, and a similar procedure will work. 

\medskip
Step 1. Obtain the group $G=PSU(5,3)$ with 

\smallskip
{\tt > G:=AtlasGroup(``U5(3)'');}

\smallskip
\noindent $G$ has order $258 190 571 520$, it has 150 conjugacy classes.  Find the two conjugacy classes of elements of order $2$, identify the larger of these two classes, which has size $444 690$.  Let $s_0$ be a representative of this conjugacy class. 

\medskip
Step 2: Let $K = C_G(s_0)$.  $K$ has five conjugacy classes of elements of order $2$, identify the largest of these, which has size $378$.  Let $s_2$ be a representative of this class. 

\medskip
Step 3: For each of the 150 conjugacy classes of $G$, find a representative $d$ of the class, set $s_1 = s_0^d$, and check if the group $\langle s_0, s_1, s_2 \rangle$ is all of $G$.  In these instances verify the intersection condition on the generating set and determine the Schl\"afli type, this verifies $G$ is a string $C$-group, with several rank $3$ representations.  

\medskip
Up to conjugacy, our implementation of Step 3 returns $36$ rank $3$ string $C$-representations of $G$.  
All of the $m_i$ occuring in their Schl\"afli types lie in the set $\{5,8,9,10,12\}$.  The Schl\"afli type $\{5,5\}$ occurs in just one case.  (For a full list of the rank $3$ Schl\"afi types for $PSU(5,3)$ along with generating permutations for each occurrence up to conjugacy, see \cite{Leemans}.)   

\medskip
Now we check the properties of irreducible characters.  We take $S=\{s_0,s_1,s_2\}$ to be the generators for the group $G$ in the string $C$-group representation of type $\{5,5\}$ generated in Step 3.  Let $H = \langle s_1,s_2 \rangle$ be the vertex stabilizer subgroup relative to this representation.  $H$ will be a dihedral group of order $10$. 

\medskip
Step 4: Generate the Atlas character table with 

\smallskip
{\tt > T:=CharacterTable(``U5(3)'');} 

\smallskip
\noindent Compute the Frobenius-Schur indicators of the irreducible characters of $G$ with 

\smallskip
{\tt > Indicator(T,2);}

\smallskip
\noindent There will be 5 irreducible characters with indicator $-1$, these are the ones with real Schur index $2$.  One of these is the irreducible character with minimal nontrivial degree $60$. 

\medskip
Step 5:  Generate the character table of $H$ and the trivial character with 

\smallskip
{\tt > T1:=CharacterTable(H);}
 
{\tt > phi:=Irr(T1)[1];}

\smallskip
\noindent Check that {\tt phi} is the trivial character of $H$. 

\medskip
Step 6: To induce $1_H$ to $G$ we need to use GAP's option to induce from table to suptable, so we obtain $(1_H)^G$ with 

\smallskip
{\tt > psi:=InducedClassFunction(T1,phi,T);} 

\medskip
Step 7: To show every nontrivial irreducible character of $G$ is polytopal, we let $\chi$ run through $Irr(G)$ and compute the multiplicities $(\chi,(1_H)^G)$ with 

\smallskip
{\tt > ScalarProduct(chi,psi);} 

\smallskip
\noindent where {\tt chi:=Irr(T)[i];} for $i=2,\dots,150$. Since $H$ is very small relative to $G$, the induced character is very close to the regular character of $G$, so it is no surprise all of these multiplicities are positive integers.  The minimal nontrivial degree character has the smallest multiplicity other than that of the trivial character of $G$, and this multiplicity is $8$. 

\bigskip
\noindent {\bf 4. Cuspidal Unipotent Characters of $PSU(n,q)$.}

\medskip
Lusztig defined cuspidal unipotent irreducible characters for groups of Lie type in \cite{Lusztig}. 

\begin{thm}[\cite{Lusztig}]
A simple group of Lie type ${}^2A_{n-1}(q)$ in odd characteristic has a unique irreducible cuspidal unipotent character if and only if $n = \frac12 s(s+1)$ for some positive integer $s$. 
\end{thm} 

Note that ${}^2A_{n-1}(q)$ is isomorphic to $PSU(n,q)$ for $q$ odd, so Lusztig's result tells us these irreducible characters exist for odd $q$ when $n$ is a triangular number.  Ohmori calculated the Schur indices of these cuspidal unipotent irreducible characters in \cite{Ohmori}.  For real Schur indices his result says: 

\begin{thm}[\cite{Ohmori}] 
The real Schur index is $2$ for the irreducible cuspidal unipotent character of $PSU(n,q)$, $q$ odd and $n=\frac12 s (s+1)$ iff $r = \lfloor \frac{n}{2} \rfloor$ is odd. 
\end{thm} 

A consequence of Ohmori's result is that for each odd $q$ there are infinitely many $PSU(n,q)$'s with an irreducible character of real Schur index $2$.  For example, this will be the case whenever $n = \frac12 s (s+1)$ with $s=4k+1$ and $k$ odd.  On the other hand, a corollary to Nuzhin's main theorem in \cite{Nuzhin} tells us 

\begin{thm}
${}^2A_{n-1}(q)$ for odd $q$ has a string $C$-group representation of rank $3$ for all $n \ge 8$.   
\end{thm}

Combining these shows there is an infinite family of string $C$-groups that have irreducible characters of real Schur index $2$.  

\begin{thm} 
If $G = PSU(n,q)$, $q$ odd, with $n = \frac12 s(s+1)$ and $r = \lfloor \frac{n}{2} \rfloor$ odd, then $G$ is a string $C$-group with an irreducible character of real Schur index $2$. 
\end{thm} 

We remark that these irreducible characters are, most likely, always going to be polytopal.


\begin{thebibliography}{99}

{\footnotesize

\bibitem{BCP97}
W.~Bosma, J.~Cannon, C.~Playoust, The {M}agma {A}lgebra {S}ystem. {I}:  the user language, {\em J. Symbolic Comput.} {\bf 24}, (1997), 235--265.
 
\bibitem{GAP} The GAP Group, GAP -- Groups, Algorithms, and Programming, Version 4.10.2; 2019. (https://www.gap-system.org) 

\bibitem{Herman-Alsairafi} A.~Herman and A.~Alsairafi, Symmetric association schemes arising from abstract regular polytopes, {\it Contrib. Discrete Math.}, {\bf 16} (1), (2021), 98-116.

\bibitem{Herman-Monson} A.~Herman and B.~Monson, On the real Schur indices associated with infinite Coxeter groups, in {\it Finite Groups 2003}, Proceedings of the Gainesville Conference on Finite Groups, edited by Chat Yin Ho et al., Walter de Gruyter, Berlin, New York, 2004, 185--194. 

\bibitem{Isaacs} I.~M.~Isaacs, Character theory of Finite Groups, Academic Press, 1976. 

\bibitem{Ladisch15} F.~Ladisch, Realizations of abstract regular polytopes from a representation theoretic view, {\it Aequationes Math.}, {\bf 90} (6), (2016), p. 1169--1193.  {\tt  arXiv:1604.07066v2 [math.MG]}

\bibitem{Leemans} D.~Leemans, Polytopes of rank 3 of $PSU(5,3)$, {\tt (http://homepages.ulb.ac.be/$\sim$dleemans/polytopes/psu53rk3.html)}

\bibitem{Lusztig} G.~Lusztig, Irreducible representations of finite classical groups, {\it Invent. Math.}, {\bf 43}, (1977), 125--175.

\bibitem{McMullen-Monson}  P.~McMullen and B.~ Monson, Realizations of regular polytopes. II. {\it Aequationes Math.}, {\bf 65} (1-2), (2003), 102--112. 

\bibitem{Nuzhin} Ya.~ N.~Nuzhin, Generating triples of involutions of Lie-type groups over a finite field of odd characteristic. I (Russian), {\it Algebra i Logika}, {\bf 36} (1), (1997), 77–96, 118; translation in {\it Algebra and Logic}, {\bf 36} (1), (1997), 46–59.

\bibitem{Ohmori} Z.~Ohmori, The Schur indices of the cuspidal unipotent characters of the finite unitary groups, 
{\it Proc. Japan Acad. Ser. A Math. Sci.}, {\bf 72} (6), (1996), 111–113. 

\bibitem{SW06} E.~Schulte and A.~I.~Weiss, Problems on Polytopes, their groups, and realizations, {\it Periodica Math. Hungarica}, {\bf 53} (2006), 231--255. 

\bibitem{AtlasRep} R.~A.~Wilson, R.~ A.~Parker, S.~Nickerson, J.~N.~Bray, and T.~Breuer, {\tt AtlasRep}, {\it A  GAP  Interface  to  the Atlas of Group Representations}, Version 1.5.0, 2019. ({\tt http://www.math.rwth-aachen.de/$\sim$Thomas.Breuer/atlasrep})

\bibitem{Atlas} R.~A.~Wilson, P.~Walsh, J.~Tripp, I.~Suleiman, R.~A.~Parker, S.~P.~Norton, S.~Nickerson, S.~Linton, J.~Bray, and R.~Abbott, ATLAS of Finite Group Representations, 2019. ({\tt http://brauer.maths.qmul.ac.uk/Atlas/v3.})

}

\end{thebibliography}
\end{document}